\documentclass[preprint,12pt]{elsarticle}
\usepackage{amsfonts}
\usepackage{amsmath}
\usepackage{amssymb}

\input{tcilatex}

\begin{document}

\begin{frontmatter}



\title{
Bayesian abstract fuzzy economies, random quasi-variational inequalities with random fuzzy mappings and random fixed point theorems
}

\author{Monica Patriche}

\address{
University of Bucharest
Faculty of Mathematics and Computer Science
    
14 Academiei Street
   
 010014 Bucharest, 
Romania
    
monica.patriche@yahoo.com }

\begin{abstract}
In this paper, we introduce an abstract fuzzy economy model with a measure space of agents 
which generalizes Patriche model (2009), we obtain a theorem of fuzzy equilibrium existence 
and we prove the existence of the solutions for two types of random quasi-variational inequalities 
with random fuzzy mappings. As a consequence, we obtain random fixed point theorems.
\end{abstract}

\begin{keyword}
Bayesian abstract fuzzy economy, \
Bayesian fuzzy equilibrium, \
incomplet information, \
random fixed point, \
random quasi-variational inequalities,\
 random fuzzy mapping.\

%% MSC codes here, in the form: \MSC code \sep code\

\end{keyword}

\end{frontmatter}

% \linenumbers


\label{}



%% References


\bibliographystyle{elsarticle-num}
\bibliography{<your-bib-database>}






% \end{thebibliography}

AMS\ Subject Classification: 58E35, 47H10, 91B50, 91A44.

\section{\textbf{\ INTRODUCTION}}

The concept of a fuzzy game (or a fuzzy abstract economy) has been
introduced in [11] and the existence of the equlibrium for 1-person fuzzy
game has been proven. The theory of fuzzy sets, initiated by Zadeh [29], has
became a good framework for obtaining results concerning fuzzy equlibrium
existence for abstract fuzzy economies. The existence of equilibrium points
of fuzzy games has been also studied in [2], [9]-[13],[17]-[19], [23].

Firstly, we introduce a fuzzy extension of Patriche's model of the Bayesian
abstract economy [16] and we prove the existence of Bayesian fuzzy
equilibrium. All economic activity in a society is made under conditions of
uncertainty (incomplet information) and the model we propose captures this
meaning. Actually, we define an abstract fuzzy economy with asymmetric
information and a measure set of agents, each of which is characterized by a
private information set, a fuzzy action (strategy) mapping, a random fuzzy
constraint one and a random fuzzy preference mapping. The fuzzy equilibrium
concept is an extension of the deterministic equilibrium. Our model
generalizes the former deterministic ones introduced by Debreu [3], Shafer
and Sonnenschein [21] or Yannelis and Prabhakar [25] or Patriche [16]. For
an overview of the results concerning the equilibrium of abstract economies,
the reader is referred to [20].

As application, we study random quasi-variational inequalities with random
fuzzy mappings. The variational inequalities were introduced in 1960s by
Fichera and Stampacchia, who studied equilibrium problems arising from
mechanics. Since then, this domain has been extensively studied and has been
found very useful in many diverse fields of pure and applied sciences, such
as mechanics, physics, optimization and control theory, operations research
and several branches of engineering sciences. Recently, the random
variational inequality and random quasi-variational inequality problems have
been introduced and studied in [7], [8], [15], [22], [28].

In this paper, we introduce an abstract fuzzy economy model with a measure
space of agents which generalizes Patriche model (2009), we prove a theorem
of fuzzy equilibrium existence and prove the existence of the solution for
two types of random quasi-variational inequalities with random fuzzy
mappings. As a consequence, we obtain random fixed point theorems.

The paper is organized as follows. In the next section, some notational and
terminological conventions are given. We also present, for the reader's
convenience, some results on Bochner integration. In Section 3, the model of
differential information abstract fuzzy economy and the main result are
presented. Section 4 contains existence theorems for solutions of random
quasi-variational inequalities with random fuzzy mappings.

\section{\textbf{\protect\smallskip NOTATION AND DEFINITION}}

Throughout this paper, we shall use the following notation:

1. $%
%TCIMACRO{\U{211d} }%
%BeginExpansion
\mathbb{R}
%EndExpansion
_{++}$ denots the set of strictly positive reals. co$D$ denotes the convex
hull of the set $D$. $\overline{co}D$ denotes the closed convex hull of the
set $D$. $2^{D}$ denotes the set of all non-empty subsets of the set $D$. If 
$D\subset Y$, where $Y$ is a topological space, cl$D$ denotes the closure of 
$D$.\smallskip

For the reader's convenience, we review a few basic definitions and results
from continuity and measurability of correspondences, Bochner integrable
functions and the integral of a correspondence.

Let $Z$ and $Y$ be sets.

\begin{definition}
The \textit{graph} of the correspondence $P:Z\rightarrow 2^{Y}$ is the set $%
G_{P}=\{(z,y)\in Z\times Y:y\in P(z)\}$.
\end{definition}

Let $Z$, $Y$ be topological spaces and $P:Z\rightarrow 2^{Y}$ be a
correspondence.

\QTP{Body Math}
1. $P$ is said to be \textit{upper semicontinuous} if for each $z\in Z$ and
each open set $V$ in $Y$ with $P(z)\subset V$, there exists an open
neighborhood $U$ of $z$ in $Z$ such that $P(y)\subset V$ for each $y\in U$.

2. $P$ is said to be \textit{lower semicontinuous} if for each $z\in Z$ and
each open set $V$ in $Y$ with $P(z)\cap V\neq \emptyset $, there exists an
open neighborhood $U$ of $z$ in $Z$ such that $P(y)\cap V\neq \emptyset $
for each $y\in U$.\medskip

\begin{lemma}
(see [27]). \textit{Let }$Z$\textit{\ and }$Y$\textit{\ be two topological
spaces and let }$D$\textit{\ be an open subset of }$Z.$\textit{\ Suppose }$%
P_{1}:Z\rightarrow 2^{Y}$\textit{\ , }$P_{2}:Z\rightarrow 2^{Y}$\textit{\
are upper.semicontinuous. such that }$P_{2}(z)\subset P_{1}(z)$\textit{\ for
all }$z\in D.$\textit{\ Then the correspondence }$P:Z\rightarrow 2^{Y}$%
\textit{\ defined by}
\end{lemma}

\begin{center}
$P\mathit{(z)=}\left\{ 
\begin{array}{c}
P_{1}(z)\text{, \ \ \ \ \ \ \ if }z\notin D\text{, } \\ 
P_{2}(z)\text{, \ \ \ \ \ \ \ \ \ \ if }z\in D%
\end{array}%
\right. $
\end{center}

\textit{is also upper semicontinuous.\medskip }

\begin{definition}
Let $Y$ be a metric space and $Y^{\prime }$ be its dual. $P:Y\rightarrow
2^{Y^{\prime }}$ is said to be\textit{\ monotone }if Re$\langle
u-v,y-x\rangle \geq 0$\textit{\ }for all $u\in P(y)$\ and $v\in P(x)$ and $%
x,y\in Y.\medskip $
\end{definition}

Let now $(\Omega ,$ $\tciFourier $, $\mu )$ be a complete, finite measure
space, and $Y$ be a topological space.

\begin{definition}
The correspondence $P:\Omega \rightarrow 2^{Y}$ is said to have a \textit{%
measurable graph} if $G_{P}\in \tciFourier \otimes \beta (Y)$, where $\beta
(Y)$ denotes the Borel $\sigma $-algebra on $Y$ and $\otimes $ denotes the
product $\sigma $-algebra.
\end{definition}

\begin{lemma}
(see\textbf{\ [}14\textbf{]).} Let $P_{n}:\Omega \rightarrow 2^{Y}$, $%
n=1,2...$be a sequence of correspondences with measurable graphs. Then the
correspondences $\cup _{n}P_{n}$, $\cap _{n}P_{n}$ and $Y\setminus P_{n}$
have measurable graphs.\medskip 
\end{lemma}

Let $(\Omega ,\tciFourier $, $\mu )$ be a measure space and $Y$ be a Banach
space.

It is known (see [14], Theorem 2, p.45) that, if $x:\Omega \rightarrow Y$ is
a $\mu $-measurable function, then $x$ is Bochner integrable if only if $%
\underset{\Omega }{\int }\Vert x(\omega )\Vert d\mu (\omega )<\infty $.

It is denoted by $L_{1}(\mu ,Y)$ the space of equivalence classes of $Y$%
-valued Bochner integrable functions $x:\Omega \rightarrow Y$ normed by $%
\parallel x\parallel =\underset{\Omega }{\int }\Vert x(\omega )\Vert d\mu
(\omega )$. Also it is known (see [4], p.50) that \textit{\ }$L_{1}(\mu ,Y)$%
\textit{\ }is a Banach space.

\begin{definition}
The correspondence $P:\Omega \rightarrow 2^{Y}$ is said to be \textit{%
integrably bounded} if there exists a map $h\in L_{1}(\mu ,\mathbb{R})$ such
that $sup\{\parallel x\parallel $ $:$ $x\in P(\omega )\}\leq h(\omega )$ $%
\mu -a.e$.
\end{definition}

We denote by $S_{P}^{1}$ the set of all selections of the correspondence $%
P:\Omega \rightarrow 2^{Y}$ that belong to the space $L_{1}$($\mu ,Y)$, i.e.

$S_{P}^{1}=\{x\in L_{1}(\mu ,Y):$ $x(\omega )\in P(\omega )$ $\mu $%
-a.e.\}.\medskip

Further, we will see the conditions under which $S_{P}^{1}$ is nonempty and
weakly compact in $L_{1}(\mu ,Y)$. Aumann measurable selection theorem (see
Appendix) and Diestel's Theorem (see Appendix) are necessary.\medskip

Let $\tciFourier (Y)$ be a collection of all fuzzy sets over $Y.$

\begin{definition}
A mapping $P:\Omega \rightarrow \tciFourier (Y)$ is called a fuzzy mapping.
If $P$ is a fuzzy mapping from $\Omega ,$ $P(\omega )$ is a fuzzy set on $Y$
and $P(\omega )(y)$ is the membership function of $y$ in $P(\omega ).$
\end{definition}

Let $A\in \tciFourier (Y),$ $\alpha \in \lbrack 0,1],$ then the set $%
(A)_{a}=\{y\in Y:A(y)\geq a\}$ is called an $a-$cut set of fuzzy set $A.$

\begin{definition}
A fuzzy mapping $P:\Omega \rightarrow \tciFourier (Y)$ is said to be
measurable if for any given $a\in \lbrack 0,1],$ $(P(\cdot ))_{a}:\Omega
\rightarrow 2^{Y}$ is a measurable set-valued mapping.
\end{definition}

\begin{definition}
We say that a fuzzy mapping $P:\Omega \rightarrow \tciFourier (Y)$ is said
to have a measurable graph if for any given $a\in \lbrack 0,1],$ the
set-valued mapping $(P(\cdot ))_{a}:\Omega \rightarrow 2^{Y}$ has a
measurable graph.
\end{definition}

\begin{definition}
A fuzzy mapping $P:\Omega \times X\rightarrow \tciFourier (Y)$ is called a
random fuzzy mapping if for any given $x\in X,$ $P(\cdot ,x):\Omega
\rightarrow \tciFourier (Y)$ is a measurable fuzzy mapping.
\end{definition}

\section{FUZZY EQUILIBRIUM EXISTENCE FOR BAYESIAN ABSTRACT FUZZY ECONOMIES}

\subsection{THE\ MODEL OF AN ABSTRACT FUZZY ECONOMY WITH A MEASURABLE SPACE
OF AGENTS}

We introduce the next model of the abstract fuzzy economy which generalizes
the classical deterministic models introduced by Debreu [3], Shafer and
Sonnenschein [21], Yannelis and Prabhakar [25] and Patriche [16].

Let $(\Omega $\textit{, }$\tciFourier ,\mu )$ be a complete finite separable
measure space, where $\Omega $ denotes the set of states of nature of the
world and the $\sigma -$algebra $\tciFourier $ denotes the set of events.
Let $Y$ be a separable Banach space, denoting the commodity or strategy
space.

Let us consider $X:T\times \Omega \rightarrow \mathcal{F}(Y)$ and $z\in
(0,1].$

Let $S_{X_{t}}^{1}=\{y(t)\in L_{1}(\mu ,Y):y(t,$\textperiodcentered $%
):\Omega \rightarrow Y$ is $\tciFourier _{t}$ $-$measurable and $y(t,\omega
)\in (X(t,\omega ))_{z}$ $\mu -a.e.\}.$ Notice that $S_{X_{t}}^{1}$ is the
set of all Bocner integrabile and $\tciFourier _{t}-$measurable selections
from the random strategy of agent $t$. In essence this is the set, out of
which agent $t$ will pick his/her optimal choices. In particular, an element 
$x_{t}$ in $S_{X_{t}}^{1}$ is called \textit{a strategy} for agent $t$. The
typical element of $S_{X_{t}}^{1}$ is denoted by $\widetilde{x}_{t}$ and
that of $X(t,\omega )$ by $x_{t}(\omega )$. Let $S_{X}^{1}=\{\widetilde{x}%
\in L_{1}(\upsilon ,L_{1}(\mu ,Y)):\widetilde{x}(t,\text{\textperiodcentered 
})\in S_{X_{t}}^{1}$ $\upsilon -a.e.\}$. An element of $S_{X}^{1}$ will be a 
\textit{joint strategy profile}.\medskip

We introduce the following model.

\begin{definition}
A \textit{Bayesian abstract fuzzy economy} (or social system) with
differential information \ and a measure space of agents $(T,\tau ,\upsilon )
$ is a set $G=\{(X,\tciFourier _{t},A,P,a,b,z),t\in T\}$, where
\end{definition}

1) $X:T\times \Omega \rightarrow \mathcal{F}(Y)$ is the random fuzzy action
(strategy) mapping, where, $X(t,\omega )\subset \mathcal{F}(Y)$ is
interpreted as the strategy set of agent $t$ of the state of nature $\omega
; $

2) $\tciFourier _{t}$ is a sub $\sigma -$algebra of $\tciFourier $ which
denotes the private information of agent $t$;

3) for each $t\in T,$ $A(t,\cdot ,\cdot ):\Omega \times S_{X}^{1}\rightarrow 
\mathcal{F}(Y)$ is the random fuzzy constraint mapping of agent $t;$

4) for each $t\in T,$ $P(t,\cdot ,\cdot ):\Omega \times S_{X}^{1}\rightarrow 
\mathcal{F}(Y)$ is the random fuzzy preference mapping of agent $t;$

5) $z\in (0,1]$ such that for all $(t,\omega ,x)\in T\times \Omega \times
S_{X}^{1},$ $(A(t,\omega ,\widetilde{x}))_{a(\widetilde{x})}\subset
(X(t,\omega ))_{z}$ and $(P(t,\omega ,\widetilde{x}))_{p(\widetilde{x}%
)}\subset X(t,\omega ))_{z};$

5) $a:S_{X}^{1}\rightarrow (0,1]$ is a random fuzzy constraint function and $%
p:S_{X}^{1}\rightarrow (0,1]$ is a random fuzzy preference function.$%
\medskip $

We propose the next concept of fuzzy equilibrium for an abstract fuzzy
economy with a measure space of agents.

\begin{definition}
\textit{A Bayesian fuzzy \ equilibrium} for $G$ is a strategy profile $%
\widetilde{x}^{\ast }\in S_{X}^{1}$ such that for $\upsilon -a.e.$,
\end{definition}

i) $\widetilde{x}^{\ast }(t,\omega )\in (A(t,\omega ,\widetilde{x}^{\ast
}))_{a(\widetilde{x}^{\ast })}$ $\mu -a.e;$

ii) $(A(t,\omega ,\widetilde{x}^{\ast }))_{a(\widetilde{x}^{\ast })}\cap
(P(t,\omega ,\widetilde{x}^{\ast }))_{p(\widetilde{x}^{\ast })}=\emptyset $ $%
\mu -a.e.\medskip $

\subsection{EXISTENCE OF THE\ BAYESIAN FUZZY\ EQUILIBRIUM}

Now, we are establishing the next equilibrium existence theorem for Bayesian
abstract fuzzy economies with a measure space of agents and with upper
semi-continuous fuzzy mappings. Our result generalizes that one we proved in
[16].

Let $(\Omega $\textit{, }$\tciFourier ,\mu )$ be a complete finite measure
separable space, where $\Omega $ denotes the set of states of nature of the
world and the $\sigma -$algebra $\tciFourier $ denotes the set of events.
Let $Y$ denote the strategy or com$\func{mod}$ity space, where $Y$ is a
separable Banach space.\medskip

\begin{theorem}
\textit{Let }$(T,\tau ,\upsilon )$\textit{\ be a measure space of agents and
be the Bayesian fuzzy abstract economy }$G=\{(X,\tciFourier
_{t},A,P,a,p,z),t\in T\}$ \textit{satisfying (A.1)-(A.5). Then, there exists
a Bayesian equilibrium for }$G$\textit{.}
\end{theorem}

\textit{A.1)}

\textit{\ \ \ \ \ (a) }$X:T\times \Omega \rightarrow \tciFourier (Y)$ 
\textit{is such that }$(t,\omega )\rightarrow (X(t,\omega ))_{z}:T\times
\Omega \rightarrow 2^{Y}$ \textit{is a non-empty convex weakly
compact-valued and integrably bounded fuzzy mapping;.}

\textit{\ \ \ \ \ (b) For each fixed }$t\in T,$\textit{\ }$(X(t,$\textit{%
\textperiodcentered }$))_{z}$\textit{\ has an }$\tciFourier _{t}-$\textit{%
measurable graph;}

\textit{\ A.2)}

\textit{\ \ \ \ \ (a) the correspondence} $(t,\omega ,\widetilde{x}%
)\rightarrow (A(t,\omega ,\widetilde{x}))_{a(\widetilde{x})}:T\times \Omega
\times S_{X}^{1}\rightarrow 2^{Y}$\textit{\ has a measurable graph;}

\textit{\ \ \ \ \ (b) For each }$(t,\omega )\in T\times \Omega ,$\textit{\ }$%
\widetilde{x}\rightarrow (A(t,\omega ,\widetilde{x}))_{a(\widetilde{x}%
)}:S_{X}^{1}\rightarrow 2^{Y}$\textit{\ is an upper semicontinuous
correspondence with closed, convex and non-empty values.}

\textit{A.3)}

\textit{\ \ \ \ (a) the correspondence} $(t,\omega ,\widetilde{x}%
)\rightarrow (P(t,\omega ,\widetilde{x}))_{p(\widetilde{x})}:T\times \Omega
\times S_{X}^{1}\rightarrow 2^{Y}$\textit{\ has a measurable graph;}

\textit{\ \ \ \ \ (b) For each }$(t,\omega )\in T\times \Omega ,$\textit{\ }$%
\widetilde{x}\rightarrow (P(t,\omega ,\widetilde{x}))_{p(\widetilde{x}%
)}:S_{X}^{1}\rightarrow 2^{Y}$\textit{\ is an upper semicontinuous
correspondence with closed, convex and non-empty values.}

\textit{A.4) The correspondence }$t\rightarrow S_{X_{t}}^{1}$\textit{\ has a
measurable graph.}

\textit{A.5)}

\textit{\ \ \ \ (a) For each }$(t,\omega )\in T\times \Omega ,$\textit{\ for
each }$\widetilde{x}\in S_{X}^{1},$\textit{\ }$\widetilde{x}(t,\omega
)\notin (A(t,\omega ,\widetilde{x}))_{a(\widetilde{x})}\cap (P(t,\omega ,%
\widetilde{x}))_{p(\widetilde{x})}$

\textit{\ \ \ \ (b) the set }$U^{(t,\omega )}=\{\widetilde{x}\in S_{X}^{1}:$%
\textit{\ }$(A(t,\omega ,\widetilde{x}))_{a(\widetilde{x})}\cap (P(t,\omega ,%
\widetilde{x}))_{p(\widetilde{x})}\neq \emptyset \}$\textit{\ is weakly open
in }$S_{X}^{1}$\textit{.}\medskip

\textit{Proof.}\textbf{\ }Define $\Phi :T\times \Omega \times
S_{X}^{1}\rightarrow 2^{Y}$ by $\Phi (t,\omega ,\widetilde{x})=(A(t,\omega ,%
\widetilde{x}))_{a(\widetilde{x})}\cap (P(t,\omega ,\widetilde{x}))_{p(%
\widetilde{x})}.$ We shall proove first that $S_{X}^{1}$ is non-empty,
convex, weakly compact.

Since $(\Omega $\textit{, }$\tciFourier ,\mu )$ is a complete finite measure
space, $Y$ is a separable Banach space and $X(t,$\textperiodcentered $%
):\Omega \rightarrow 2^{Y}$ has measurable graph, by Aumann's selection
theorem (see Appendix) it follows that there exists a function $f(t,$%
\textperiodcentered $):\Omega \rightarrow Y$ such that $f(t,\omega )\in
X(t,\omega )$ $\mu -a.e$. Since $X(t,$\textperiodcentered $)$ is integrably
bounded, we have that $f(t,$\textperiodcentered $)\in L_{1}(\mu ,Y)$, hence $%
S_{X_{t}}$ is non-empty and $S_{X}=\tprod\limits_{t\in T}S_{X_{t}}$ is
non-empty. $S_{X_{t}}^{1}$ is convex and $S_{X}^{1}$ is also convex. Since $%
X(t,$\textperiodcentered $):\Omega \rightarrow 2^{Y}$ is integrably bounded
and has convex weakly compact values, by Diestel's Theorem (see Appendix) it
follows that $S_{X_{t}}^{1}$ is a weakly compact subset of $L_{1}(\mu ,Y)$
and so is $S_{X}^{1}.$ We have that $S_{X}^{1}$ is a metrizable set as being
a weakly compact subset of the separable Banach space $L_{1}(\upsilon
,L_{1}(\mu ,Y))$ (Dunford-Schwartz [1958, p 434]).

Let denote $A^{\ast }:T\times \Omega \times S_{X}^{1}\rightarrow 2^{Y},$ $%
A^{\ast }((t,\omega ,\widetilde{x}))=(A(t,\omega ,\widetilde{x}))_{a(%
\widetilde{x})}.$

\QTP{Body Math}
Since all the values of the correspondence $A^{\ast }$ are contained in the
compact set $X($\textperiodcentered $,$\textperiodcentered $)$ and $A^{\ast
} $ is closed and convex valued (hence weakly closed), it follows that $%
A^{\ast }$ is weakly compact valued.

\QTP{Body Math}
Then $\Phi $ is convex valued and for each $(t,\omega )\in T\times \Omega ,$ 
$\Phi (t,\omega ,$\textit{\textperiodcentered }$)$ is upper semicontinuous.
We have that $\Phi ($\textperiodcentered ,\textperiodcentered , $\widetilde{x%
})$ has a measurable graph for each $\widetilde{x}\in S_{X}^{1}$. Let $%
U=\{(t,\omega ,\widetilde{x})\in \Omega \times S_{X}^{1}:\Phi (t,\omega ,%
\widetilde{x})\neq \emptyset \}$. For each $\widetilde{x}\in S_{X}^{1},$\ let%
\textit{\ }$U^{x}=\{(t,\omega )\in T\times \Omega :\Phi (t,\omega ,%
\widetilde{x})\neq \emptyset \}$ and for each $\omega \in \Omega ,$\ let $%
U^{(t,\omega )}=\{\widetilde{x}\in S_{X}^{1}:\Phi (t,\omega ,\widetilde{x}%
)\neq \emptyset \}$. Define $G:T\times \Omega \times S_{X}^{1}\rightarrow
2^{Y}$ by $G(t,\omega ,\widetilde{x})=\left\{ 
\begin{array}{c}
\Phi (t,\omega ,\widetilde{x})\text{ if }(t,\omega ,\widetilde{x})\in U \\ 
(A(t,\omega ,\widetilde{x}))_{a(\widetilde{x})}\text{ if }(t,\omega ,%
\widetilde{x})\notin U.%
\end{array}%
\right. $

\QTP{Body Math}
For each $\widetilde{x}\in S_{X}^{1}$,the correspondence $G($%
\textperiodcentered ,\textperiodcentered $,\widetilde{x})$ has a measurable
graph.

\QTP{Body Math}
By assumption A5) (b), the set $U^{(t,\omega )}=\{\widetilde{x}\in
S_{X}^{1}: $\ $(A(t,\omega ,\widetilde{x}))_{a(\widetilde{x})}\cap
(P(t,\omega ,\widetilde{x}))_{p(\widetilde{x})}\neq \emptyset \}$ is weakly
open in $S_{X}^{1}$. For each $(t,\omega )\in T\times \Omega ,$ $G(t,\omega
, $\textit{\textperiodcentered }$):S_{X}^{1}\rightarrow 2^{Y}$ is upper
semi-continuous.

Let $V$ be weakly open in $S_{X}^{1}$ and $\widetilde{x}\in S_{X}^{1}$.

$W=\{\widetilde{x}\in S_{X}^{1}:G(t,\omega ,\widetilde{x})\subset V\}=$

\ \ \ $=\{\widetilde{x}\in U^{(t,\omega )}:G(t,\omega ,\widetilde{x})\subset
V\}\cup \{\widetilde{x}\in S_{X}^{1}\setminus U^{(t,\omega )}:G(t,\omega ,%
\widetilde{x})\subset V\}$

\ \ \ $=\{\widetilde{x}\in U^{(t,\omega )}:\Phi (t,\omega ,\widetilde{x}%
)\subset V\}\cup \{\widetilde{x}\in S_{X}^{1}:G(t,\omega ,\widetilde{x}%
)\subset V\}.$

$W$ is an open set, because $U^{(t,\omega )}$ is open, $\Phi (t,\omega ,$%
\textperiodcentered $)$ is a upper semicontinuous map on $U^{(t,\omega )}$
and the set $\left\{ \widetilde{x}\in S_{X}^{1}:\text{ }G(t,\omega ,%
\widetilde{x})\subset V\right\} $ is open since $G(t,\omega ,$%
\textperiodcentered $)$ is u.s.c. Moreover, $G$ is convex and non-empty
valued.

\QTP{Body Math}
Define the correspondence $\varphi :T\times S_{X}^{1}\rightarrow
2^{L_{1}(\mu ,Y)}$ by $\varphi (t,\widetilde{x})=\{\widetilde{y}(t)\in
L_{1}(\mu ,Y):\widetilde{y}(t,\omega )\in G(t,\omega ,\widetilde{x})$ $\mu
-a.e.\}\cap S_{X_{t}}^{1}.$

\QTP{Body Math}
By the measurability lifting theorem (Theorem 16, see Appendix), the
correspondence $t\rightarrow \{\widetilde{y}(t)\in L_{1}(\mu ,Y):\widetilde{y%
}(t,\omega )\in G(t,\omega ,\widetilde{x})$ $\mu -a.e.\}$ has a measurable
graph and so does $t\rightarrow S_{X_{t}}^{1}$ by (A.4). Thus, for each
fixed $\widetilde{x}\in S_{X}^{1}$, $\varphi ($\textperiodcentered $,%
\widetilde{x})$ has a measurable graph.

\QTP{Body Math}
$\func{Si}$nce for each fixed $\widetilde{x}\in S_{X}^{1}$, $G($\textit{%
\textperiodcentered }$,$\textit{\textperiodcentered }$,\widetilde{x})$ has a
measurable graph and it is nonempty valued, then by the Aumann measurable
selection theorem (see Appendix), it admits a measurable slection and we can
conclude that $\varphi $\ is nonempty valued. It folows by the u.s.c.
litfing theorem that for each fixed t, $\varphi (t,$\textperiodcentered $)$
is weakly u.s.c.

\QTP{Body Math}
Define $G^{^{\prime }}:S_{X}^{1}\rightarrow 2^{S_{X}^{1}}$, by $G^{^{\prime
}}(\widetilde{x})=\{\widetilde{y}\in S_{X}^{1}:\widetilde{y}(t)\in \varphi
(t,\widetilde{x})$ $\nu -a.e.\}.$

\QTP{Body Math}
Another application of the u.s.c. lifting theorem enables us to conclude
that $G^{^{\prime }}$ is a weakly u.s.c. correspondence which is obviously
convex valued (since $\varphi $ is convex valued) and also nonempty valued
(recall once more the Aumann measurable selection theorem and notice that
the set $S_{X}^{1}$ is metrizable).

$G^{^{\prime }}$ is an upper semicontinuous correspondence and has also
non-empty convex closed values.

By Fan-Glicksberg's fixed-point theorem in [6], there exists $\widetilde{x}%
^{\ast }\in S_{X}^{1}$ such that $\widetilde{x}^{\ast }\in G^{^{\prime }}(%
\widetilde{x}^{\ast })$. It follows that $\widetilde{x}^{\ast }\in S_{X}^{1}$
and $\widetilde{x}^{\ast }(t)\in \varphi (t,\widetilde{x}^{\ast })$ $\nu
-a.e.$ Thus, we have that $\widetilde{x}^{\ast }\in S_{X}^{1}$ and $%
\widetilde{x}^{\ast }(t,\omega )\in G(t,\omega ,\widetilde{x}^{\ast })$ $\mu
-a.e.,$ $\nu -a.e.$

By assumption A.4.a, it follows that $\widetilde{x}^{\ast }(t,\omega )\notin
(A(t,\omega ,\widetilde{x}^{\ast }))_{a(\widetilde{x}^{\ast })}\cap
(P(t,\omega ,\widetilde{x}^{\ast }))_{p(\widetilde{x}^{\ast })}$, then we
have that $\widetilde{x}^{\ast }\notin U$ and $\widetilde{x}^{\ast
}(t,\omega )\in (A(t,\omega ,\widetilde{x}^{\ast }))_{a(\widetilde{x}^{\ast
})}.$

Therefore, for $\upsilon -a.e.:$

i) $\widetilde{x}^{\ast }(t,\omega )\in (A(t,\omega ,\widetilde{x}^{\ast
}))_{a(\widetilde{x}^{\ast })}$ $\mu -a.e;$

ii) $(A(t,\omega ,\widetilde{x}^{\ast }))_{a(\widetilde{x}^{\ast })}\cap
(P(t,\omega ,\widetilde{x}^{\ast }))_{p(\widetilde{x}^{\ast })}=\emptyset $ $%
\mu -a.e.\Box \medskip $

\section{RANDOM\ QUASI-VARIATIONAL INEQUALITIES}

In this section, we are establishing new random quasi-variational
inequalities with random fuzzy mappings and random fixed point theorems. The
proofs rely on the theorem of Bayesian fuzzy equilibrium existence for the
abstract fuzzy economy with a measure space of agents.

This is our first theorem.

\begin{theorem}
\textit{Let }$(T,\tau ,\upsilon )$\textit{\ be a measure space, }$(\Omega $%
\textit{, }$\tciFourier ,\mu )$ be a complete finite separable measure
space, and let $Y$ be a separable Banach space.\textit{\ Suppose that the
following conditions are satisfied:}
\end{theorem}

\textit{A.1)}

\textit{\ \ \ \ \ (a) }$X:T\times \Omega \rightarrow F(Y)$ and $z\in (0,1]$ 
\textit{are such that} $(t,\omega )\rightarrow (X(t,\omega ))_{z}:T\times
\Omega \rightarrow 2^{Y}$ \textit{is a non-empty convex weakly
compact-valued and integrably bounded correspondence;}

\textit{\ \ \ \ \ (b) For each fixed }$t\in T,$\textit{\ }$\tciFourier _{t}$ 
\textit{is a sub }$\sigma -$\textit{algebra of }$\tciFourier $\textit{\ and} 
$(X(t,$\textit{\textperiodcentered }$))_{z}$\textit{\ has an }$\tciFourier
_{t}-$\textit{measurable graph;}

\textit{\ A.2)}

\textit{\ \ \ \ \ }$A:T\times \Omega \times S_{X}^{1}\rightarrow \mathcal{F}%
(Y)$\textit{\ and }$a:S_{X}^{1}\rightarrow (0,1]$ \textit{are such that }

\ \ \ \ \textit{(a)} \textit{the correspondence} $(t,\omega ,\widetilde{x}%
)\rightarrow (A(t,\omega ,\widetilde{x}))_{a(\widetilde{x})}:T\times \Omega
\times S_{X}^{1}\rightarrow 2^{Y}$\textit{\ has a measurable graph;}

\textit{\ \ \ \ \ (b) For each }$(t,\omega )\in T\times \Omega ,$\textit{\ }$%
\widetilde{x}\rightarrow (A(t,\omega ,\widetilde{x}))_{a(\widetilde{x}%
)}:S_{X}^{1}\rightarrow 2^{Y}$\textit{\ is an upper semicontinuous
correspondence with closed, convex and non-empty values.}

\textit{\ \ \ \ A.3) The correspondence }$t\rightarrow S_{X_{t}}^{1}$\textit{%
\ has a measurable graph;}

\textit{A.4)}

\ \ \ $\psi :T\times \Omega \times S_{X}^{1}\times Y\rightarrow \mathbb{R}%
\cup \{-\infty ,+\infty \}$\textit{\ and }$p:S_{X}^{1}\rightarrow (0,1]$ 
\textit{are} \textit{such that:}

\ \ \ (\textit{a) }$(t,\widetilde{x})\rightarrow \{y\in Y:\psi (t,\omega ,%
\widetilde{x},y)>0\}$\textit{\ is upper semicontinuous with compact values
on }$T\times S_{X}^{1}$\textit{\ for each fixed }$\omega \in \Omega ;$

\ \ \ (\textit{b) }$\widetilde{x}(t,\omega )\notin $cl$\{y\in Y:\psi
(t,\omega ,\widetilde{x},y)>0\}$\textit{\ for each fixed }$(t,\omega ,%
\widetilde{x})\in T\times \Omega \times S_{X}^{1};$

\ \ \ (\textit{c) for each} $(t,\omega ,\widetilde{x})\in T\times \Omega
\times S_{X}^{1},$ $\psi (t,\omega ,\widetilde{x},\cdot )$ \textit{is
concave;}

\ \ \ (\textit{d) for each }$(t,\omega )\in T\times \Omega ,$\textit{\ }$\{%
\widetilde{x}\in S_{X}^{1}\mathit{\ }:\alpha (t,\omega ,\widetilde{x})>0\}$%
\textit{\ is weakly open in} $S_{X}^{1}\mathit{\ },$\textit{\ where }$\alpha
:T\times \Omega \times S_{X}^{1}\rightarrow \mathbb{R}$\textit{\ is defined
by }$\alpha (t,\omega ,\widetilde{x})=\sup_{y\in (A(t,\omega ,\widetilde{x}%
))_{a(\widetilde{x})}}\psi (t,\omega ,\widetilde{x},y)$\textit{\ for each }$%
(t,\omega ,\widetilde{x})\in T\times \Omega \times S_{X}^{1};$

\ \ \ (\textit{d) }$\{(t,\omega ,\widetilde{x}):\psi (t,\omega ,\widetilde{x}%
,y)>0\}\in \tau \otimes \mathcal{F}_{t}\otimes B(S_{X}^{1})$\textit{.}

\textit{Then, there exists }$\widetilde{x}^{\ast }\in S_{X}^{1}$\textit{\
such that }$\upsilon -a.e.$,

i) $\widetilde{x}^{\ast }(t,\omega )\in (A(t,\omega ,\widetilde{x}^{\ast
}))_{a(\widetilde{x}^{\ast })}$ $\mu -a.e;$

\textit{ii) sup}$_{y\in (A(t,\omega ,\widetilde{x}^{\ast }))_{a(\widetilde{x}%
^{\ast })}}\psi (t,\omega ,\widetilde{x}^{\ast },y)\leq 0$\textit{\ }$\mu
-a.e.\medskip $

\textit{Proof.} Let $P:T\times \Omega \times S_{X}^{1}\rightarrow \mathcal{F}%
(Y)$ and $p:S_{X}^{1}\rightarrow (0,1]$ such that $(P(t,\omega ,\widetilde{x}%
))_{p(\widetilde{x})}=\{y\in Y:\psi (t,\omega ,\widetilde{x},y)>0\}$ for
each $(t,\omega ,\widetilde{x})\in T\times \Omega \times S_{X}^{1}.$

We shall show that the abstract fuzzy economy $G=\{(X,\tciFourier
_{t},A,P,a,p,z),t\in T\}$ satisfies all hypotheses of Theorem 4.

Suppose $(t,\omega )\in T\times \Omega .$

According to A4 a), we have that\textit{\ }for each\textit{\ }$(t,\omega
)\in T\times \Omega ,$\textit{\ }$\widetilde{x}\rightarrow (P(t,\omega ,%
\widetilde{x}))_{p(\widetilde{x})}:S_{X}^{1}\rightarrow 2^{Y}$\textit{\ }is
an upper semicontinuous correspondence with closed, convex and non-empty
values and according to A4 b), $\widetilde{x}(t,\omega )\not\in (P(t,\omega ,%
\widetilde{x}))_{p(\widetilde{x})}$ for each $\widetilde{x}\in S_{X}^{1}.$
Assumption A4 c) implies that $\widetilde{x}\rightarrow (P(t,\omega ,%
\widetilde{x}))_{p(\widetilde{x})}:S_{X}^{1}\rightarrow 2^{Y}$ has convex
values.

According to the definition of $\alpha ,$ we note that $\{\widetilde{x}\in
S_{X}^{1}:(A(t,\omega ,\widetilde{x}^{\ast }))_{a(\widetilde{x}^{\ast
})}\cap $

\noindent $\cap (P(t,\omega ,\widetilde{x}^{\ast }))_{p(\widetilde{x}^{\ast
})}\neq \emptyset \}=\{\widetilde{x}\in S_{X}^{1}:\alpha (t,\omega ,%
\widetilde{x})>0\}$ so that $\{\widetilde{x}\in S_{X}^{1}:(A(t,\omega ,%
\widetilde{x}^{\ast }))_{a(\widetilde{x}^{\ast })}\cap (P(t,\omega ,%
\widetilde{x}^{\ast }))_{p(\widetilde{x}^{\ast })}\neq \emptyset \}$ is
weakly open in $S_{X}^{1}$ by A4d).

According to A2 a) and A4 e), it follows that the correspondences\textit{\ }$%
(t,\omega ,\widetilde{x})\rightarrow (A(t,\omega ,\widetilde{x}^{\ast }))_{a(%
\widetilde{x}^{\ast })}$ and $(t,\omega ,\widetilde{x})\rightarrow
(P(t,\omega ,\widetilde{x}^{\ast }))_{p(\widetilde{x}^{\ast })}$\textit{\ }%
have measurable graphs$.$

Thus the abstract fuzzy economy $G=\{(X,\tciFourier _{t},A,P,a,p,z),t\in T\}$
satisfies all hypotheses of Theorem 3. Therefore, there exists $\widetilde{x}%
^{\ast }\in S_{X}^{1}$ such that $\upsilon -a.e.:$

$x^{\ast }(t,\omega )\in (A(t,\omega ,\widetilde{x}^{\ast }))_{a(\widetilde{x%
}^{\ast })}$ $\mu -a.e$ and

$(A(t,\omega ,\widetilde{x}^{\ast }))_{a(\widetilde{x}^{\ast })}\cap
(P(t,\omega ,\widetilde{x}^{\ast }))_{p(\widetilde{x}^{\ast })}=\phi $ $\mu
-a.e;$

that is, there exists $\widetilde{x}^{\ast }\in S_{X}^{1}$ such that $%
\upsilon -a.e.:$

i) $x^{\ast }(t,\omega )\in (A(t,\omega ,\widetilde{x}^{\ast }))_{a(%
\widetilde{x}^{\ast })}$ $\mu -a.e;$

ii) sup$_{y\in (A(t,\omega ,\widetilde{x}^{\ast }))_{a(\widetilde{x}^{\ast
})}}\psi (t,\omega ,\widetilde{x}^{\ast },y)\leq 0$ $\mu -a.e.\medskip $

If \TEXTsymbol{\vert}T\TEXTsymbol{\vert}=1, we obtain the following
corollary.

\begin{corollary}
\textit{Let }$(\Omega $\textit{, }$\tciFourier ,\mu )$ be a complete finite
separable measure space, and let $Y$ be a separable Banach space.\textit{\
Suppose that the following conditions are satisfied:}
\end{corollary}

\textit{A.1)}

\textit{\ \ \ \ \ (a) }$X:\Omega \rightarrow F(Y)$ and $z\in (0,1]$ \textit{%
are such that} $\omega \rightarrow (X(\omega ))_{z}:T\times \Omega
\rightarrow 2^{Y}$ \textit{is a non-empty convex weakly compact-valued and
integrably bounded correspondence;}

\textit{\ \ \ \ \ (b) }$(X($\textit{\textperiodcentered }$))_{z}$\textit{\
has an }$\tciFourier -$\textit{measurable graph;}

\textit{\ A.2)}

\textit{\ \ \ \ \ }$A:\Omega \times S_{X}^{1}\rightarrow \mathcal{F}(Y)$%
\textit{\ and }$a:S_{X}^{1}\rightarrow (0,1]$ \textit{are such that }

\ \ \ \ \textit{(a)} \textit{the correspondence} $(\omega ,\widetilde{x}%
)\rightarrow (A(\omega ,\widetilde{x}))_{a(\widetilde{x})}:\Omega \times
S_{X}^{1}\rightarrow 2^{Y}$\textit{\ has a measurable graph;}

\textit{\ \ \ \ \ (b) For each }$\omega \in \Omega ,$\textit{\ }$\widetilde{x%
}\rightarrow (A(\omega ,\widetilde{x}))_{a(\widetilde{x})}:S_{X}^{1}%
\rightarrow 2^{Y}$\textit{\ is an upper semicontinuous correspondence with
closed, convex and non-empty values.}

\textit{\ A.3)}

\ \ \ $\psi :\Omega \times S_{X}^{1}\times Y\rightarrow \mathbb{R}\cup
\{-\infty ,+\infty \}$\textit{\ and }$p:S_{X}^{1}\rightarrow (0,1]$ \textit{%
are} \textit{such that:}

\ \ \ (\textit{a) }$\widetilde{x}\rightarrow \{y\in Y:\psi (\omega ,%
\widetilde{x},y)>0\}$\textit{\ is upper semicontinuous with compact values
on }$S_{X}^{1}$\textit{\ for each fixed }$\omega \in \Omega ;$

\ \ \ (\textit{b) }$\widetilde{x}(\omega )\notin $cl$\{y\in Y:\psi (\omega ,%
\widetilde{x},y)>0\}$\textit{\ for each fixed }$(\omega ,\widetilde{x})\in
\Omega \times S_{X}^{1};$

\ \ \ (\textit{c) for each} $(\omega ,\widetilde{x})\in \Omega \times
S_{X}^{1},$ $\psi (\omega ,\widetilde{x},\cdot )$ \textit{is concave;}

\ \ \ (\textit{d) for each }$\omega \in \Omega ,$\textit{\ }$\{\widetilde{x}%
\in S_{X}^{1}\mathit{\ }:\alpha (\omega ,\widetilde{x})>0\}$\textit{\ is
weakly open in} $S_{X}^{1}\mathit{\ },$\textit{\ where }$\alpha :\Omega
\times S_{X}^{1}\rightarrow \mathbb{R}$\textit{\ is defined by }$\alpha
(\omega ,\widetilde{x})=\sup_{y\in (A(\omega ,\widetilde{x}))_{a(\widetilde{x%
})}}\psi (\omega ,\widetilde{x},y)$\textit{\ for each }$(\omega ,\widetilde{x%
})\in \Omega \times S_{X}^{1};$

\ \ \ (\textit{d) }$\{(\omega ,\widetilde{x}):\psi (\omega ,\widetilde{x}%
,y)>0\}\in \mathcal{F}_{t}\otimes B(S_{X}^{1})$\textit{.}

\textit{Then, there exists }$\widetilde{x}^{\ast }\in S_{X}^{1}$\textit{\
such that:}

i) $\widetilde{x}^{\ast }(t,\omega )\in (A(t,\omega ,\widetilde{x}^{\ast
}))_{a(\widetilde{x}^{\ast })}$ $\mu -a.e;$

\textit{ii) sup}$_{y\in (A(\omega ,\widetilde{x}^{\ast }))_{a(\widetilde{x}%
^{\ast })}}\psi (\omega ,\widetilde{x}^{\ast },y)\leq 0$\textit{\ }$\mu
-a.e.\medskip $

As a consequence of Theorem 4, we prove the following random generalized
quasi-variational inequality with random fuzzy mappings. Our result
generalizes Theorem 5.1 in [22].

\begin{theorem}
\textit{Let }$(T,\tau ,\upsilon )$\textit{\ be a measure space, }$(\Omega $%
\textit{, }$\tciFourier ,\mu )$ be a complete finite separable measure
space, and let $Y$ be a separable Banach space.\textit{\ Suppose that the
following conditions are satisfied:}
\end{theorem}

\textit{A.1)}

\textit{\ \ \ \ \ (a) }$X:T\times \Omega \rightarrow F(Y)$ and $z\in (0,1]$ 
\textit{are such that} $(t,\omega )\rightarrow (X(t,\omega ))_{z}:T\times
\Omega \rightarrow 2^{Y}$ \textit{is a non-empty convex weakly
compact-valued and integrably bounded correspondence;}

\textit{\ \ \ \ \ (b) For each fixed }$t\in T,$\textit{\ }$\tciFourier _{t}$
is a sub $\sigma -$algebra of $\tciFourier $ and $(X(t,$\textit{%
\textperiodcentered }$))_{z}$\textit{\ has an }$\tciFourier _{t}-$\textit{%
measurable graph;}

\textit{\ A.2)}

\textit{\ \ \ \ \ }$A:T\times \Omega \times S_{X}^{1}\rightarrow \mathcal{F}%
(Y)$\textit{\ and }$a:S_{X}^{1}\rightarrow (0,1]$ \textit{are such that }

\ \ \ \ \textit{(a)} \textit{the correspondence} $(t,\omega ,\widetilde{x}%
)\rightarrow (A(t,\omega ,\widetilde{x}))_{a(\widetilde{x})}:T\times \Omega
\times S_{X}^{1}\rightarrow 2^{Y}$\textit{\ has a measurable graph;}

\textit{\ \ \ \ \ (b) For each }$(t,\omega )\in T\times \Omega ,$\textit{\ }$%
\widetilde{x}\rightarrow (A(t,\omega ,\widetilde{x}))_{a(\widetilde{x}%
)}:S_{X}^{1}\rightarrow 2^{Y}$\textit{\ is an upper semicontinuous
correspondence with closed, convex and non-empty values.}

\textit{A.3) The correspondence }$t\rightarrow S_{X_{t}}^{1}$\textit{\ has a
measurable graph;}

\textit{A.4)}

\ \ \ $G:T\times \Omega \times Y\rightarrow \mathcal{F}(Y^{\prime })$ and $%
g:Y\rightarrow (0,1]$ \textit{are such that:}

\ \ \ (\textit{a) for each }$(t,\omega )\in T\times \Omega ,$\textit{\ }$%
y\rightarrow (G(t,\omega ,y))_{g(y)}:Y\rightarrow 2^{Y^{\prime }}$ \textit{%
is monotone\ with non-empty values}$;$

\ \ \ (\textit{b)for each }$(t,\omega )\in T\times \Omega ,$\textit{\ }$%
y\rightarrow (G(t,\omega ,y))_{g(y)}:L\cap Y\rightarrow 2^{Y^{\prime }}$%
\textit{\ is lower semicontinuous from the relative topology of }$Y$\textit{%
\ into the weak}$^{\ast }-$\textit{topology }$\sigma (Y^{\prime },Y)$ 
\textit{of }$Y^{\prime }$\textit{\ for each one-dimensional flat }$L\subset
Y;$

\textit{A.5)}

\ \ \ (\textit{a) }$f:T\times \Omega \times S_{X}^{1}\times Y\rightarrow 
\mathbb{R}\cup \{\infty ,-\infty \}$\textit{\ is such that }$\widetilde{x}%
\rightarrow f(t,\omega ,\widetilde{x,}y)$\textit{\ is lower semicontinuous
on }$S_{X}^{1}$\textit{\ for each fixed }$(t,\omega ,y)\in T\times \Omega
\times Y,$\textit{\ }$f(t,\omega ,\widetilde{x},\widetilde{x}(t,\omega ))=0$%
\textit{\ for each }$(t,\omega ,\widetilde{x})\in T\times \Omega \times
S_{X}^{1}$\textit{\ and }$y\rightarrow f(t,\omega ,\widetilde{x,}y)$\textit{%
\ is concave on }$Y$ \textit{for each fixed }$(t,\omega ,x)\in T\times
\Omega \times S_{X}^{1}$\textit{;}

\ \ \ (\textit{b) for each fixed }$(t,\omega )\in T\times \Omega ,$\textit{\
the set}

\textit{\ }$\{\widetilde{x}\in S_{X}^{1}:\sup_{y\in (A(t,\omega ,\widetilde{x%
}))_{a(\widetilde{x})}}[\sup_{u\in (G(t,\omega ,y))_{g(y)}}$\textit{Re}$%
\langle u,\widetilde{x}-y\rangle +f(t,\omega ,\widetilde{x},y)]>0\}$\textit{%
\ is weakly open in }$S_{X}^{1}$

\ \ \ (c) $\{(t,\omega ,\widetilde{x}):\sup_{u\in (G(t,\omega ,y))_{g(y)}}%
\mathit{Re}\langle u,\widetilde{x}-y\rangle +f(t,\omega ,\widetilde{x}%
,y)>0\}\in \tau \otimes \mathcal{F}_{t}\otimes B(S_{X}^{1})$\textit{.}

\textit{Then, there exists }$\widetilde{x}^{\ast }\in S_{X}^{1}$\textit{\
such that} $\upsilon -a.e.$:

i) $\widetilde{x}^{\ast }(t,\omega )\in (A(t,\omega ,\widetilde{x}^{\ast
}))_{a(\widetilde{x}^{\ast })}$ $\mu -a.e;$

ii) \textit{sup}$_{u\in (G(t,\omega ,y))_{g(y)}}Re\langle u,x^{\ast }(\omega
)-y\rangle +f(t,\omega ,xy)]\leq 0$\textit{\ for all }$y\in (A(t,\omega ,%
\widetilde{x}^{\ast }))_{a(\widetilde{x}^{\ast })}$\textit{\ }$\mu
-a.e.\medskip $

\textit{Proof.} Let us define $\psi :T\times \Omega \times S_{X}^{1}\times
Y\rightarrow \mathbb{R}\cup \{-\infty ,+\infty \}$ by

$\psi (t,\omega ,\widetilde{x},y)=\sup_{u\in (G(t,\omega ,y))_{g(y)}}$%
\textit{Re}$\langle u,\widetilde{x}-y\rangle +f(t,\omega ,\widetilde{x},y)$
for each $(t,\omega ,\widetilde{x},y)\in T\times \Omega \times
S_{X}^{1}\times Y.$

According to assumption A5 a), $\widetilde{x}\rightarrow f(t,\omega ,%
\widetilde{x},y)$ is lower semicontinuous on $S_{X}^{1}$ for each fixed $%
(t,\omega ,y)\in T\times \Omega \times Y$ and $f(t,\omega ,\widetilde{x},%
\widetilde{x}(t,\omega ))=0$ for each $(t,\omega ,\widetilde{x})\in T\times
\Omega \times S_{X}^{1}$ implies that $\widetilde{x}(t,\omega )\notin \{y\in
Y:\psi (t,\omega ,\widetilde{x},y)>0\}$\textit{\ }for each fixed\textit{\ }$%
(t,\omega ,\widetilde{x})\in T\times \Omega \times S_{X}^{1}.$

We also have that for each $(t,\omega ,\widetilde{x})\in T\times \Omega
\times S_{X}^{1},$ $\psi (t,\omega ,\widetilde{x},\cdot )$ is concave$.$
This fact is a consequence of assumption A4 a).

All the hypotheses of Theorem 4 are satisfied. According to Theorem 4, there
exists $\widetilde{x}^{\ast }\in S_{X}^{1}$ such that $\upsilon -a.e.$

$\widetilde{x}^{\ast }(t,\omega )\in (A(t,\omega ,\widetilde{x}^{\ast }))_{a(%
\widetilde{x}^{\ast })}$ $\mu -a.e.,$

and

(1) \ \ sup$_{y\in (A(t,\omega ,\widetilde{x}^{\ast }))_{a(\widetilde{x}%
^{\ast })}}\sup_{(G(t,\omega ,y))_{g(y)}}[\mathit{Re}\langle u,\widetilde{x}%
^{\ast }(t,\omega )-y\rangle +f(t,\omega ,\widetilde{x}^{\ast },y)]\leq 0$ $%
\mu -a.e.$

Finally, we will prove that $\upsilon -a.e.,$

sup$_{y\in A(t,\omega ,\widetilde{x}^{\ast }))_{a(\widetilde{x}^{\ast
})}}\sup_{u\in (G(t,\omega ,\widetilde{x}^{\ast }))_{g(x^{\ast })}}[\mathit{%
Re}\langle u,\widetilde{x}^{\ast }(\omega )-y\rangle +f(t,\omega ,\widetilde{%
x}^{\ast },y)]\leq 0$ $\mu -a.e.$

In order to do that, let us consider the fixed points $t\in T$ and $\omega
\in \Omega .$

Let $y\in v$, $\lambda \in \lbrack 0,1]$ and $z_{\lambda }(t,\omega
):=\lambda y+(1-\lambda )\widetilde{x}^{\ast }(t,\omega ).$ According to
assumption A2 b), $z_{\lambda }(t,\omega )\in (A(t,\omega ,\widetilde{x}%
^{\ast }))_{a(\widetilde{x}^{\ast })}.$

According to (1), we have $\sup_{u\in (G(t,\omega ,z_{\lambda }(t,\omega
)))_{g(z_{\lambda }(t,\omega ))}}[\mathit{Re}\langle u,\widetilde{x}^{\ast
}(t,\omega )-z_{\lambda }(t,\omega )\rangle +f(t,\omega ,\widetilde{x}^{\ast
},z_{\lambda }(t,\omega ))]\leq 0$ for each $\lambda \in \lbrack 0,1]$.

According to assumption A5 a), $f(t,\omega ,\widetilde{x}^{\ast },\widetilde{%
x}^{\ast }(t,\omega ))=0$\textit{.} For each $y_{1},y_{2}\in Y$\ and for
each $\lambda \in \lbrack 0,1],$\ we also have that $f(t,\omega ,\widetilde{x%
}^{\ast },\lambda y_{1}+(1-\lambda )y_{2})\geq \lambda f(t,\omega ,%
\widetilde{x}^{\ast },y_{1})+(1-\lambda )f(t,\omega ,\widetilde{x}^{\ast
},y_{2}).$

Therefore, for each $\lambda \in \lbrack 0,1]$, we have that

$t\{\sup_{u\in (G(t,\omega ,z_{\lambda }(t,\omega )))_{g(z_{\lambda
}(t,\omega ))}}[\mathit{Re}\langle u,\widetilde{x}^{\ast }(t,\omega
)-y\rangle +f(t,\omega ,\widetilde{x}^{\ast },y)]\}\leq $

$\sup_{u\in (G(t,\omega ,z_{\lambda }(t,\omega )))_{g(z_{\lambda }(t,\omega
))}}t[\mathit{Re}\langle u,\widetilde{x}^{\ast }(t,\omega _{i})-y)\rangle
+f(t,\omega ,\widetilde{x}^{\ast },z_{\lambda }(t,\omega ))]=$

$\sup_{u\in (G(t,\omega ,z_{\lambda }(t,\omega )))_{g(z_{\lambda }(t,\omega
))}}[\mathit{Re}\langle u,\widetilde{x}^{\ast }(t,\omega )-z_{\lambda
}(t,\omega )\rangle +f(t,\omega ,\widetilde{x}^{\ast },z_{\lambda }(t,\omega
))]\leq 0.$

It follows that for each $\lambda \in \lbrack 0,1],$

(2) $\sup_{u\in (G(t,\omega ,z_{\lambda }(t,\omega )))_{g(z_{\lambda
}(t,\omega ))}}[\mathit{Re}\langle u,\widetilde{x}^{\ast }(t,\omega
)-y\rangle +f(t,\omega ,\widetilde{x}^{\ast },y)]\leq 0.$

Now, we are using the lower semicontinuity of $F(t,\omega ,\cdot )$ in order
to show the conclusion. For each $z_{0}\in (G(t,\omega ,\widetilde{x}^{\ast
}(t,\omega )))_{g(\widetilde{x}^{\ast }(t,\omega ))}$ and $e>0$ let us
consider $U_{z_{0}}^{t},$ the neighborhood of $z_{0}$ in the topology $%
\sigma (Y^{\prime },Y),$ defined by $U_{z_{0}}^{t}:=\{z\in Y^{\prime }:|%
\func{Re}\langle z_{0}-z,\widetilde{x}^{\ast }(t,\omega )-y\rangle |<e\}.$
As $y\rightarrow (G(t,\omega ,y))_{g(y)}:L\cap Y\rightarrow 2^{Y^{\prime }}$
is lower semicontinuous, where $L=\{z_{\lambda }(t,\omega ):\lambda \in
\lbrack 0,1]\}$ and $U_{z_{0}}^{t}\cap (G(t,\omega ,\widetilde{x}^{\ast
}(t,\omega )))_{g(\widetilde{x}^{\ast }(t,\omega ))}\neq \emptyset ,$ there
exists a non-empty neighborhood $N(\widetilde{x}^{\ast }(t,\omega ))$ of $%
\widetilde{x}^{\ast }(t,\omega )$ in $L$ such that for each $z\in N(%
\widetilde{x}^{\ast }(t,\omega )),$ we have that $U_{z_{0}}^{t}\cap
F(t,\omega ,z)\neq \emptyset .$ Then there exists $\delta \in (0,1],$ $t\in
(0,\delta )$ and $u\in (G(t,\omega ,z_{\lambda }(t,\omega )))_{g(z_{\lambda
}(t,\omega ))}\cap U_{z_{0}}^{t}\neq \emptyset $ such that $\mathit{Re}%
\langle z_{0}-u,\widetilde{x}^{\ast }(t,\omega )-y\rangle <e.$ Therefore, $%
\mathit{Re}\langle z_{0},\widetilde{x}^{\ast }(t,\omega )-y\rangle <\mathit{%
Re}\langle u_{i},\widetilde{x}^{\ast }(t,\omega )-y\rangle +e.$

It follows that

$\mathit{Re}\langle z_{0},\widetilde{x}^{\ast }(t,\omega )-y\rangle
+f(t,\omega ,\widetilde{x}^{\ast },y)<\mathit{Re}\langle u,\widetilde{x}%
^{\ast }(t,\omega )-y\rangle +f(t,\omega ,\widetilde{x}^{\ast },y)+e<e.$

The last inequality comes from (2). Since $e>0$ and $z_{0}\in (G(t,\omega ,%
\widetilde{x}^{\ast }(t,\omega )))_{g(\widetilde{x}^{\ast }(t,\omega ))}$
have been chosen arbitrarily, the next relation holds:

$\mathit{Re}\langle z_{0},\widetilde{x}^{\ast }(t,\omega )-y\rangle
+f(t,\omega ,\widetilde{x}^{\ast },y)<0.$

Hence, $\mu -a.e.,$ we have that $\sup_{u\in (G(t,\omega ,\widetilde{x}%
^{\ast }(t,\omega )))_{g(\widetilde{x}^{\ast }(t,\omega ))}}$[$\mathit{Re}%
\langle z_{0},\widetilde{x}^{\ast }(t,\omega )-y\rangle +f(t,\omega ,%
\widetilde{x}^{\ast },y)]\leq 0$ for every $y\in (A(t,\omega ,\widetilde{x}%
^{\ast }))_{a(\widetilde{x}^{\ast })}$.$\medskip $

If \TEXTsymbol{\vert}T\TEXTsymbol{\vert}=1, we can make abstraction of T and
we obtain the following corollary.

\begin{theorem}
\textit{Let }$(\Omega $\textit{, }$\tciFourier ,\mu )$ be a complete finite
separable measure space, and let $Y$ be a separable Banach space.\textit{\
Suppose that the following conditions are satisfied:}
\end{theorem}

\textit{A.1)}

\textit{\ \ \ \ \ (a) }$X:\Omega \rightarrow F(Y)$ and $z\in (0,1]$ \textit{%
are such that} $\omega \rightarrow (X(\omega ))_{z}:\Omega \rightarrow 2^{Y}$
\textit{is a non-empty convex weakly compact-valued and integrably bounded
correspondence;}

\textit{\ \ \ \ \ (b) }$(X(,$\textit{\textperiodcentered }$))_{z}$\textit{\
has an }$\tciFourier -$\textit{measurable graph;}

\textit{\ A.2)}

\textit{\ \ \ \ \ }$A:\Omega \times S_{X}^{1}\rightarrow \mathcal{F}(Y)$%
\textit{\ and }$a:S_{X}^{1}\rightarrow (0,1]$ \textit{are such that }

\ \ \ \ \textit{(a)} \textit{the correspondence} $(\omega ,\widetilde{x}%
)\rightarrow (A(\omega ,\widetilde{x}))_{a(\widetilde{x})}:T\times \Omega
\times S_{X}^{1}\rightarrow 2^{Y}$\textit{\ has a measurable graph;}

\textit{\ \ \ \ \ (b) For each }$\omega \in \Omega ,$\textit{\ }$\widetilde{x%
}\rightarrow (A(\omega ,\widetilde{x}))_{a(\widetilde{x})}:S_{X}^{1}%
\rightarrow 2^{Y}$\textit{\ is an upper semicontinuous correspondence with
closed, convex and non-empty values;}

\textit{A.3)}

\ \ \ $G:\Omega \times Y\rightarrow \mathcal{F}(Y^{\prime })$ and $%
g:Y\rightarrow (0,1]$ \textit{are such that:}

\ \ \ (\textit{a) for each }$\omega \in \Omega ,$\textit{\ }$y\rightarrow
(G(\omega ,y))_{g(y)}:Y\rightarrow 2^{Y^{\prime }}$ \textit{is monotone\
with non-empty values}$;$

\ \ \ (\textit{b)for each }$\omega \in \Omega ,$\textit{\ }$y\rightarrow
(G(\omega ,y))_{g(y)}:L\cap Y\rightarrow 2^{Y^{\prime }}$\textit{\ is lower
semicontinuous from the relative topology of }$Y$\textit{\ into the weak}$%
^{\ast }-$\textit{topology }$\sigma (Y^{\prime },Y)$ \textit{of }$Y^{\prime
} $\textit{\ for each one-dimensional flat }$L\subset Y;$

\textit{A.5)}

\ \ \ (\textit{a) }$f:\Omega \times S_{X}^{1}\times Y\rightarrow \mathbb{R}%
\cup \{\infty ,-\infty \}$\textit{\ is such that }$\widetilde{x}\rightarrow
f(\omega ,\widetilde{x,}y)$\textit{\ is lower semicontinuous on }$S_{X}^{1}$%
\textit{\ for each fixed }$(\omega ,y)\in \Omega \times Y,$\textit{\ }$%
f(\omega ,\widetilde{x},\widetilde{x}(t,\omega ))=0$\textit{\ for each }$%
(\omega ,\widetilde{x})\in \Omega \times S_{X}^{1}$\textit{\ and }$%
y\rightarrow f(\omega ,\widetilde{x,}y)$\textit{\ is concave on }$Y$ \textit{%
for each fixed }$(\omega ,x)\in \Omega \times S_{X}^{1}$\textit{;}

\ \ \ (\textit{b) for each fixed }$\omega \in \Omega ,$\textit{\ the set}

\textit{\ }$\{\widetilde{x}\in S_{X}^{1}:\sup_{y\in (A(\omega ,\widetilde{x}%
))_{a(\widetilde{x})}}[\sup_{u\in (G(\omega ,y))_{g(y)}}$\textit{Re}$\langle
u,\widetilde{x}-y\rangle +f(t,\omega ,\widetilde{x},y)]>0\}$\textit{\ is
weakly open in }$S_{X}^{1}$

\ \ \ (c) $\{(\omega ,\widetilde{x}):\sup_{u\in (G(\omega ,y))_{g(y)}}%
\mathit{Re}\langle u,\widetilde{x}-y\rangle +f(t,\omega ,\widetilde{x}%
,y)>0\}\in \mathcal{F}_{t}\otimes B(S_{X}^{1})$\textit{.}

\textit{Then, there exists }$\widetilde{x}^{\ast }\in S_{X}^{1}$\textit{\
such that}:

i) $\widetilde{x}^{\ast }(\omega )\in (A(\omega ,\widetilde{x}^{\ast }))_{a(%
\widetilde{x}^{\ast })}$ $\mu -a.e;$

ii) \textit{sup}$_{u\in (G(\omega ,y))_{g(y)}}Re\langle u,x^{\ast }(\omega
)-y\rangle +f(t,\omega ,xy)]\leq 0$\textit{\ for all }$y\in (A(\omega ,%
\widetilde{x}^{\ast }))_{a(\widetilde{x}^{\ast })}$\textit{\ }$\mu
-a.e.\medskip $

We obtain the following random fixed point theorem with fuzzy mapping as a
particular case of Theorem 6$.$

\begin{theorem}
\textit{Let }$(T,\tau ,\upsilon )$\textit{\ be a measure space, }$(\Omega $%
\textit{, }$\tciFourier ,\mu )$ be a complete finite separable measure
space, and let $Y$ be a separable Banach space.\textit{\ Suppose that the
following conditions are satisfied:}
\end{theorem}

\textit{A.1)}

\textit{\ \ \ \ \ (a) }$X:T\times \Omega \rightarrow F(Y)$ and $z\in (0,1]$ 
\textit{are such that} $(t,\omega )\rightarrow (X(t,\omega ))_{z}:T\times
\Omega \rightarrow 2^{Y}$ \textit{is a non-empty convex weakly
compact-valued and integrably bounded correspondence;}

\textit{\ \ \ \ \ (b) For each fixed }$t\in T,$\textit{\ }$\tciFourier _{t}$
is a sub $\sigma -$algebra of $\tciFourier $ and $(X(t,$\textit{%
\textperiodcentered }$))_{z}$\textit{\ has an }$\tciFourier _{t}-$\textit{%
measurable graph;}

\textit{\ A.2)}

\textit{\ \ \ \ \ }$A:T\times \Omega \times S_{X}^{1}\rightarrow \mathcal{F}%
(Y)$\textit{\ and }$a:S_{X}^{1}\rightarrow (0,1]$ \textit{are such that }

\ \ \ \ \textit{(a)} \textit{the correspondence} $(t,\omega ,\widetilde{x}%
)\rightarrow (A(t,\omega ,\widetilde{x}))_{a(\widetilde{x})}:T\times \Omega
\times S_{X}^{1}\rightarrow 2^{Y}$\textit{\ has a measurable graph;}

\textit{\ \ \ \ \ (b) For each }$(t,\omega )\in T\times \Omega ,$\textit{\ }$%
\widetilde{x}\rightarrow (A(t,\omega ,\widetilde{x}))_{a(\widetilde{x}%
)}:S_{X}^{1}\rightarrow 2^{Y}$\textit{\ is an upper semicontinuous
correspondence with closed, convex and non-empty values.}

\textit{A.3) The correspondence }$t\rightarrow S_{X_{t}}^{1}$\textit{\ has a
measurable graph;}

\textit{Then, there exists }$\widetilde{x}^{\ast }\in S_{X}^{1}$\textit{\
such that} $\upsilon -a.e.$, $\widetilde{x}^{\ast }(t,\omega )\in
(A(t,\omega ,\widetilde{x}^{\ast }))_{a(\widetilde{x}^{\ast })}$ $\mu
-a.e.\medskip $

If \TEXTsymbol{\vert}T\TEXTsymbol{\vert}=1, we obtain the following result.

\begin{theorem}
\textit{Let }$(\Omega $\textit{, }$\tciFourier ,\mu )$ be a complete finite
separable measure space, and let $Y$ be a separable Banach space.\textit{\
Suppose that the following conditions are satisfied:}
\end{theorem}

\textit{A.1)}

\textit{\ \ \ \ \ (a) }$X:\Omega \rightarrow F(Y)$ and $z\in (0,1]$ \textit{%
are such that} $\omega \rightarrow (X(\omega ))_{z}:T\times \Omega
\rightarrow 2^{Y}$ \textit{is a non-empty convex weakly compact-valued and
integrably bounded correspondence;}

\textit{\ \ \ \ \ (b) }$(X($\textit{\textperiodcentered }$))_{z}$\textit{\
has an }$\tciFourier -$\textit{measurable graph;}

\textit{\ A.2)}

\textit{\ \ \ \ \ }$A:\Omega \times S_{X}^{1}\rightarrow \mathcal{F}(Y)$%
\textit{\ and }$a:S_{X}^{1}\rightarrow (0,1]$ \textit{are such that }

\ \ \ \ \textit{(a)} \textit{the correspondence} $(\omega ,\widetilde{x}%
)\rightarrow (A(\omega ,\widetilde{x}))_{a(\widetilde{x})}:\Omega \times
S_{X}^{1}\rightarrow 2^{Y}$\textit{\ has a measurable graph;}

\textit{\ \ \ \ \ (b) For each }$\omega \in \Omega ,$\textit{\ }$\widetilde{x%
}\rightarrow (A(\omega ,\widetilde{x}))_{a(\widetilde{x})}:S_{X}^{1}%
\rightarrow 2^{Y}$\textit{\ is an upper semicontinuous correspondence with
closed, convex and non-empty values.}

\textit{Then, there exists }$\widetilde{x}^{\ast }\in S_{X}^{1}$\textit{\
such that} $\widetilde{x}^{\ast }(\omega )\in (A(\omega ,\widetilde{x}))_{a(%
\widetilde{x})}$ $\mu -a.e.\medskip $

\section{APPENDIX}

The results below have been used in the proof of our theorems. For more
details and further references see the paper quoted.

\begin{theorem}
\textit{(Projection theorem).} \textit{Let }$(\Omega ,$\textit{\ }$%
\tciFourier $\textit{, }$\mu )$\textit{\ be a complete, finite measure
space, and }$Y$\textit{\ be a complete separable metric space. If }$H$%
\textit{\ belongs to }$\tciFourier \otimes $\textit{\ss }$(Y)$\textit{, its
projection Proj}$_{\Omega }(H)$\textit{\ belongs to }$\tciFourier .$
\end{theorem}

\begin{theorem}
\textit{(Aumann measurable selection theorem [24])}. \textit{Let }$(\Omega $%
\textit{, }$\tciFourier ,\mu )$\textit{\ be a complete finite measure space, 
}$Y$\textit{\ be a complete, separable metric space and }$T:\Omega
\rightarrow 2^{Y}$\textit{\ be a nonempty valued correspondence with a
measurable graph, i.e., }$G_{T}\in \tciFourier \otimes \beta (Y).$\textit{\
Then there is a measurable function }$f:\Omega \rightarrow Y$\textit{\ such
that }$f(\omega )\in T(\omega )$\textit{\ }$\mu -a.e.\medskip $
\end{theorem}

\begin{theorem}
\textit{(Diestel's Theorem [24, Theorem 3.1).\ Let }$(\Omega $\textit{, }$%
\tciFourier ,\mu )$\textit{\ be a complete finite measure space, }$X$\textit{%
\ be a separable Banach space and }$T:\Omega \rightarrow 2^{Y}$\textit{\ be
an integrably bounded, convex, weakly compact and nonempty valued
correspondence. Then }$S_{T}=\{x\in L_{1}(\mu ,Y):$ $x(\omega )\in T(\omega )
$ $\mu $-a.e.\}\textit{\ is weakly compact in }$L_{1}(\mu ,Y).$\textit{%
\smallskip }
\end{theorem}

\begin{theorem}
\textbf{(}\textit{Carath\'{e}odory-type selection theorem} [14]\textbf{)}. 
\textit{Let }$(\Omega $\textit{,}$\tciFourier ,\mu )$\textit{\ be a complete
measure space, }$Z$\textit{\ be a complete separable metric space and }$Y$%
\textit{\ a separable Banach space. Let }$X:\Omega \rightarrow 2^{Y}$\textit{%
\ be a correspondence with a measurable graph, i.e., }$G_{X}\in \tciFourier
\otimes $\textit{\ss }$(Y)$\textit{\ and let }$T:\Omega \times Z\rightarrow
2^{Y}$\textit{\ be a convex valued correspondence (possibly empty) with a
meaurable graph, i.e., }$G_{T}\in \tciFourier \otimes $\textit{\ss }$%
(Z)\otimes $\textit{\ss }$(Y)$\textit{\ where }\text{\text{\ss }}$(Y)$%
\textit{\ and \ss }$(Z)$\textit{\ are the Borel }$\sigma -$\textit{algebras
of }$Y$\textit{\ and }$Z$\textit{, respectively}$.$
\end{theorem}

\textit{Suppose that:}

\textit{(1) for each }$\omega \in \Omega $\textit{, }$T(\omega ,x)\subset
X(\omega )$\textit{\ for all }$x\in Z.$

\textit{(2) for each }$\omega \in \Omega $\textit{, }$T(\omega ,$\textit{%
\textperiodcentered }$)$\textit{\ has open lower sections in Z, i.e., for
each }$\omega \in \Omega $\textit{\ and }$y\in Y$\textit{, }$T^{-1}(\omega
,y)=\{x\in Z:y\in T(\omega ,x)\}$\textit{\ is open in Z.}

\textit{(3) for each }$(\omega ,x)\in \Omega \times Z,$\textit{\ if }$%
T(\omega ,x)\neq \emptyset $\textit{, then }$T(\omega ,x)$\textit{\ has a
non-empty interior in }$X(\omega )$\textit{.}

\textit{Let }$U=\{(\omega ,x)\in \Omega \times Z:T(\omega ,x)\neq \emptyset
\}$\textit{\ and for each }$x\in Z$\textit{, }$U^{x}=\{\omega \in \Omega
:(\omega ,x)\in U\}$\textit{\ and for each }$\omega \in \Omega $\textit{, }$%
U^{\omega }=\{x\in Z:(\omega ,x)\in U\}.$\textit{\ Then for each }$x\in Z,$%
\textit{\ }$U^{x}$\textit{\ is a measurable set in }$\Omega $\textit{\ and
there exists a Caratheodory-type selection from }$T_{\mid U}$, \textit{\
i.e., there exists a function }$f:U\rightarrow Y$\textit{\ such that }$%
f(\omega ,x)\in T(\omega ,x)$\textit{\ for all }$(\omega ,x)\in U$\textit{\
, for each }$x\in Z,$\textit{\ }$f($\textit{\textperiodcentered }$,x)$%
\textit{\ is measurable on }$U^{x}$\textit{\ and for each }$\omega \in
\Omega ,$\textit{\ }$f(\omega ,$\textit{\textperiodcentered }$)$\textit{\ is
continuous on }$U^{\omega }.$\textit{\ Moreover, }$f($\textit{%
\textperiodcentered }$,$\textit{\textperiodcentered }$)$\textit{\ is jointly
measurable.\medskip }

\begin{theorem}
\textit{(U. s. c. Lifting Theorem. [24]).} \textit{Let }$Y$\textit{\ be a
separable space, }$(\Omega $\textit{, }$\tciFourier ,\mu )$\textit{\ be a
complete finite measure space and }$X:\Omega \rightarrow 2^{Y}$\textit{\ be
an integrably bounded, nonempty, convex valued correspondence such that for
all }$\omega \in \Omega ,$\textit{\ }$X(\omega )$\textit{\ is a weakly
compact, convex subset of }$Y$\textit{. Denote by }$S_{X}$\textit{\ the set }%
$\{x\in L_{1}(\mu ,Y):x(\omega )\in X(\omega )$\textit{\ }$\mu -a.e.\}.$%
\textit{\ Let }$T:\Omega \times S_{X}\rightarrow 2^{Y}$\textit{\ be a
nonempty, closed, convex valued correspondence such that }$T(\omega
,x)\subset X(\omega )$\textit{\ for all }$(\omega ,x)\in \Omega \times
S_{X}^{1}.$\textit{\ Assume that for each fixed }$x\in S_{X},$\textit{\ }$T($%
\textit{\textperiodcentered }$,x)$\textit{\ has a measurable graph and that
for each fixed }$\omega \in \Omega ,$\textit{\ }$T(\omega ,$\textit{%
\textperiodcentered }$):S_{X}\rightarrow 2^{Y}$\textit{\ is u.s.c. in the
sense that the set }$\{x\in S_{X}:T(\omega ,x)\subset V\}$\textit{\ is
weakly open in }$S_{X}$\textit{\ for every norm open subset }$V$\textit{\ of 
}$Y$\textit{. Define the correspondence }$\Phi :S_{X}\rightarrow 2^{S_{X}}$%
\textit{\ by}
\end{theorem}

$\Phi (x)=\{y\in S_{X}:y(\omega )\in T(\omega ,x)$\textit{\ }$\mu -a.e.\}.$

\textit{Then }$\Phi $\textit{\ is weakly u.s.c., i.e., the set }$\{x\in
S_{X}:\Phi (x)\subset V\}$\textit{\ is weakly open in }$S_{X}$\textit{\ for
every weakly open subset }$V$\textit{\ of }$S_{X}$\textit{.}

\begin{theorem}
(Measurability Lifting Theorem) [1]. \textit{Let }$Y$\textit{\ and }$E$%
\textit{\ be separable Banach spaces, and }$(T,\tau ,\upsilon )$\textit{\
and }$(\Omega ,\tciFourier \mu )$\textit{\ be finite complete separable
measure spaces. Let }$\gamma :T\times \Omega \times E\rightarrow 2^{Y}$%
\textit{\ be a nonempty valued correspondence. Suppose that for each }$y\in
E,$\textit{\ }$\gamma ($\textit{\textperiodcentered ,\textperiodcentered ,}$%
y)$\textit{\ has a measurable graph. Define the correspondence }$\psi
:\Omega \times E\rightarrow 2^{L_{1}(\mu ,Y)}$\textit{\ by}
\end{theorem}

$\psi (t,y)=\{x(t)\in L_{1}(\mu ,Y):x(t,\omega )\in \gamma (t,\omega ,y)$%
\textit{\ }$\mu -a.e.\}.$

\textit{Then for each }$y\in E,\psi ($\textit{\textperiodcentered }$,y)$%
\textit{\ has a measurable graph.}

\end{document}